# A set of formulas for primes

by
Simon Plouffe
April 4, 2022

Abstract

In 1947, W. H. Mills published a paper describing a formula that gives primes : if A = 1.3063778838630806904686144926… then $\lfloor A^{3^n} \rfloor$ is always prime, here $\lfloor x \rfloor$ is the integral part of x. Later in 1951, E. M. Wright published another formula, if $g_0 = \alpha = 1.9287800\ldots$ and $g_{n+1} = 2^{g_n}$ then

$$\lfloor g_n \rfloor = \left\lfloor 2^{\cdots 2^{2^{\alpha}}} \right\rfloor \text{ is always prime.}$$

The primes are uniquely determined by $\alpha$ , the prime sequence is 3, 13, 16381, …
The growth rate of these functions is very high since the fourth term of Wright formula is a 4932 digit prime and the 8'th prime of Mills formula is a 762 digit prime.
A new set of formulas is presented here, giving an arbitrary number of primes minimizing the growth rate. The first one is : if $a_0 = 43.8046877158\ldots$ and $a_{n+1} = a_n^{\frac{5}{4}}$, then if $S(n)$ is the rounded values of $a_n$, $S(n) = 113$, 367, 1607, 10177, 102217, 1827697, 67201679, 6084503671, …. Other exponents can also give primes like 11/10, or 101/100. If $a_0$ is well chosen then it is conjectured that any exponent $> 1$ can also give an arbitrary series of primes. When the exponent is 3/2 it is conjectured that all the primes are within a series of trees. The method for obtaining the formulas is explained. Five formulas are presented and all results are empirical.

Résumé

En 1947, W. H. Mills publiait un article montrant une formule qui peut donner un nombre arbitraire de nombres premiers. Si A = 1.3063778838630806904686144926… alors $\lfloor A^{3^n} \rfloor$ donne une suite arbitraire de nombres tous premiers. , ici $\lfloor x \rfloor$ est le plancher de x. Plus tard en 1951, E. M. Wright en proposait une autre, si $g_0 = \alpha = 1.9287800\ldots$ et $g_{n+1} = 2^{g_n}$ alors

$$\lfloor g_n \rfloor = \left\lfloor 2^{\cdots 2^{2^{\alpha}}} \right\rfloor \text{ est toujours premier.}$$

Les premiers consécutifs sont uniquement représentés par $\alpha$. La suite de premiers est 3, 13, 16381,…
Le taux de croissance de ces 2 fonctions est assez élevé puisque le 4ème terme de la suite de Wright a 4932 chiffres décimaux. La croissance de celle de Mills est moins élevée, le 8ème terme a quand même une taille de 762 chiffres. Une série de formules est présentée ici qui minimise le taux de croissance et qui possède les mêmes propriétés de fournir une suite de premiers de longueur arbitraire. Si $a_0 = 43.8046877158\ldots$ et $a_{n+1} = a_n^{\frac{5}{4}}$ alors la suite $S(n) = \{ a_n \}$ : l'arrondi de $a_n$, est une suite de premiers de longueur arbitraire. Ici l'exposant $^{5}/_{4}$ peut être abaissé à $^{11}/_{10}$, ou même $^{101}/_{100}$. Si $a_0$ est bien choisi il est conjecturé que l'exposant peut être aussi près de 1 que l'on veut. Un autre modèle est basé sur $\{c^n\}$, où { } est l'arrondi avec c = 31622.77671… et n > 1 permet d'obtenir 387 nombres premiers d'affilée. Cinq formules sont présentées et tous les résultats sont empiriques.

# Introduction

The first type of prime formula to consider is for example, given $a_0$ a real constant $> 0$ and $a_{n+1} = a(n)10$, if $a_0 = 7.3327334517988679\ldots$ then the sequence $73, 733, 7333, 73327, 733273, \ldots$ is a sequence of primes but fails for obvious reasons after a few terms. If the base is changed to any other fixed size base, taking into account that the average gap between primes is increasing then eventually the process fails to give any more primes.

If we choose a function that grows faster like $n^n$, we get better results. The best start constant found is $c = 0.265588372943143390897129453665466129438 9\ldots$ giving 19 primes.
But fails at 23 (beginning at n = 3). Here $a_n = \lfloor cn^n \rfloor$.

7
67
829
12391
218723
4455833
102894377
2655883729
75775462379
2368012611049
80440106764817
2951219812933057
116299525867995629
4899240744635092571
219705395187452015923
10449948501874965563651
525445257345556693801913
27848959374722952425334841
1551723179991864497606172809

Again, for the same reasons mentioned earlier, the process fails to go further, no better example was found. The method used is a homemade Monte-Carlo method that uses Simulated Annealing (principe du recuit simulé in French).

If a function grows too slowly, eventually the average gap between primes increases and the process ceases to give any more primes. The next step was to consider formulas like Mills or Wright. The question was then : is there a way to get a useful formula that grows just enough to produce primes ?

If we consider the recurrence $a_{n+1} = {a_n}^2 - a_n + 1$ that arises in the context of Sylvester sequence. The Sylvester sequence is A000058 of the OEIS catalogue and begins with $a_0 = 2$, like this : (2, 3, 7, 43, 1807, 3263443, 10650056950807, ), that sequence has the property that

$$1 = \frac{1}{2} + \frac{1}{3} + \frac{1}{7} + \frac{1}{43} + \frac{1}{1807} + \frac{1}{3263443} + \cdots$$

The natural extension that comes next is: can we choose $a_0$ so that $a_n$ will always produce primes ? The answer is yes, when $a(0) = 1.6181418093242092 \ldots$ and by using the $\lfloor x \rfloor$ function we get,

$$a(n) = 2, 3, 7, 43, 1811, 3277913, 10744710357637, \ldots$$

The sequence and formula are interesting for one reason the growth rate is quite smaller than the one of Mills and Wright.

## A Formula for primes

What if we choose the exponent to be as small as possible? The problem with that last one is that it is stilll growing too fast, $a(14) = 9.838 \ldots \times 10^{1667}$. The size of primes doubles in length at each step.

The simplest found was $a_0 = 43.80468771580293481\ldots$ and using the { x } , rounding to the nearest integer, we get

$$a_{n+1} = a_n^{\frac{5}{4}}$$

Now, what if we carefully choose $a_0$ so that the exponent is smaller, would it work ? Let's try with 11/10 and start with a larger number.

$$a_{n+1} = a_n^{\frac{11}{10}}$$

For example when $S_0 = 100000000000000000000000000000049.31221074776345 \ldots$ and the exponent beeing 11/10 then we get the primes :

100000000000000000000000000000049
158489319246111348520210137339236753
524807460249772597364312157022725894401
3908408957924020300919472370957356345933709
70990461585528724931289825118059422005340095813
3438111840350699188044461057631015443312900908952333
489724690004200094265557071425023036671550364178496540501
...

If we want a smaller starting value then $a_0$ has to be bigger, I could get a series of primes when $a_0 = 10^{64} + 57 + \varepsilon$ , where $0 < \varepsilon < 0.5$ chosen at random. In this case the exponent is

$$a_{n+1} = a_n^{\frac{21}{20}}$$

If we choose $a_0 = 10^{600} + 543 + \varepsilon$ then we get our formula to be.

$$a_{n+1} = a_n^{\frac{101}{100}}$$

If a(0) = 2.03823915478206876746349086260954825144862477844317361... and the exponent 3/2 then the sequence of primes is :

3
5
11
37
223
3331
192271
84308429
774116799347
681098209317971743
562101323304225290104514179
1332667822014585978282511662572214575909
1538448162271607869601834587431948506238982765193425993274489

The natural question that comes next is : Can we generate all the primes with one single exponent ? Here is the tree graph of primes with the exponent 3/2.

## Other ways of generating primes

Three other methods are presented. The first is based on a very high degree polynomial equation. We begin by the known Euler polynomial$n^2 + n + 41$ that produces 40 primes in a row when n goes from 0 to 39. In 2010, François Dress and Bernard Landreau found a $6^{th}$ degree polynomial that gives 56 primes values. It was obtained after several months of computation. The polynomial is

$$\frac{n^6}{72} - \frac{5n^5}{24} + \frac{n^6}{72} - \frac{1493n^4}{72} + \frac{1493n^3}{72} + \frac{100471n^2}{18} - \frac{11971n}{6} - 57347\,.$$

All other attempts failed to go further. If we take another approach to it there is more.

If c is a real positive constant then

$$[cn^k + 1/2]$$

Where [ . ] denotes the integer part function. For a fixed k much larger than 2 it is possible to calculate c using the simulated annealing and Monte-Carlo method used before. For example, if $k = 64$ then by using the formula

$$p = nextprime(n^k)$$

$$c = \frac{p}{n^k}$$

Then the value of the new c is calculated by rounding properly

$$p = nextprime\left(\left[cn^k + \frac{1}{2}\right]\right)$$

$$c_{k+1} = \frac{p}{(n+1)^k}$$

The hypothesis is : for a large k, c is constant and will produce an arbitrary number of primes. When k = 64, the method quickly finds 8 primes values but if k = 599 then 70 primes are produced. The value of c found is :

```
1. 0000000000000000000000000000000000000000000000000000000000000\
000000000000000000000000000000000000000000000000000000000000000\
000000000000000000000000000000000000000000000000000001556808232\
856463140060524522424280772934429200613394089802373785244587288\
947728630209290312111238140972072983869588571830870498647900351\
444778386216146076636515297094962921999030102707899655875083040\
665145525741443694679398437806398341541385226365101056503219529\
413038357761603719148231623498218255233850813147568998744093447\
755692241415245289626420200959624751354414736085310934015664680\
524107510074349656822480162844410900000000000000000000006345174\
532890027939294089368385587928671551372984708501441585620418908\
704346095361979911978459546006902805287612053866464609018035122\
943293118288150470557741929174121124290158771911527759861726269\
584165754859588606741441445142853262133802425477033893100601759\
975564856010884639189352574717896440249667645304974867391425986\
281906375105044471660347105158250104363991834652263421461752834\
975705269597512556117680247632410969440028501070105083651713662\
403142442324953222225461614957943694687384105246348009815340760\
743692542284884699931895728286412542979603098044884198867262644\
327255465732046886424762403984621514353794473387936916290029911\
549331529954509625933626714214956674379391784430822602413846681\
316909794084400273637943424856690207465817148161052758912779700\
994744869116506018623229585002806564436851770494944385807359855\
293872240395791047188062073418411994344296873461271952410865026\
038448423723192039177383879070574428090707165141606873722573750\
361413628118939617030810117859967828092642486835827304995753842\
763548295039042203138474167878091109846423833753023399379895 96\
168347847772485744109449547853544023508997274222527142776574434\
316734525917951119963956001777673781864527502261570125603162371\
258224293816247417631732394744100802942101015408972089551902455\
081699086722777589571816095705739686859606854254953133788523318\
655544571643416019528317460445413625249324735829
```

When k = 1001 I find this value of c that produces 104 prime values.

```
1. 000000000000000000000000000000000000000000000000000000000000000000000000000000000000000\
00000000000000000000000000000000000000000000000000000000000000000000000000000000000000000\
00000000000000000000000000000000000000000000000000000000000000000000000000000000000000000\
0000000000000000000000000000000000125990588497934548663662088537715317898307345260180583 3\
23931435872455532981729710651349444158213948818545323640593772754269090815382794373303605\
31032253065209164543992867389533988645378869293054521535540437640264924717218337701861582\
63863610200230076170685773794205724815861643935173361519983613116882690284144248101920208\
06510182210643573796316380910863539606120370866893232496051485086497467990763578114679066\
43511345081362364422769127194021649235153662464913703294618822166887983870315588441106214\
81340962305982173042950831271359740413862293006103829962131522304901811854369079285560801\
56127395594055635262316335676598235139936600519471752446533884127775136800820321897034833\
6925216314960153661963670000000000000000000000000000000000000002856922463537929582277227375\
34573209459155555553507901662476263229096688827287105130336980209309890445915932025320437\
91865264980534443344524609828999110761933292589504030787982875567498653229908488978554566\
67343197183608335096725513262211244819036081998055736478108166949284397076878265689251779\
50165437169403417121093904025741600510426629665273037772624863935742023884842513100551563\
13386198931648636209922198705909308260842342176288588385011256962600297234767059688581215\
42128399100716160760517497248823653252822696639553588275394359123882969533179056493207692\
39058072976876703918551210565097947101808802571780934485397548675096907719003422433055454\
36432892896196532241017067958091280226589158192434780052667595497166072344482563094303132\
07749969618831846158232622183434850148402756911211398935776546073125564968263058674825326\
78385994696717675249981709965958000987443947212244644763379515406938514672346131121294611\
46321844369295493249482895087729753260849786693705728092698851465492392173730706849962492\
31216596166252982070773148562282049797435559272711648811046626031928698550228861837334923\
47445282121848836355413332570898852324742813132203189317097096931337090249364473376188254\
80333372330490460254399100556222506183678193972395584173595952332514796665024284845488485\
64456888046546242546151999688740068991295519427373610371260392861045349960918614305187606\
01199946505452669414820340092566861924748261672447488843829277513346089163004507571313765\
06792446848015950986087309825879656484491208770711299508536138483156405005289125879215048\
```

262580884594130104731092685277684005452284638098913830988830807268502167793938868506533390\
0282953922794033452518766212590

Which is a record since the polynomial $[\, cn^{1001} + {}^{1}/_{2}]$ will produce 104 primes when n = 2..105.

The second formula is based on an observation about primes in arithmetic progression. The idea is to change the base asumption. Instead of an arithmetic progression, I use a geometric progression. If $f$ is defined as :

$$f(n) = [c^n + 1/2]$$

Then when c = 2.553854696… , 7 primes are generated : 3, 7, 17, 43, 109, 277 and 709. But what if c >> 1 ? Again by using the previous method I can produce this table of values for c

| Constant c such that $[c^n + 1/2]$ is prime | Values | Number of primes produced |
| --- | --- | --- |
| 2.553854696… | 3, 7, 17, 43, 109, 277, 709 | 7 |
| 2027.1671684764912194343956 | n=1..97 | 97 |
| 577.181936975247888… | n=1..22 | 22 |
| 593.46526943871… | n=2..48 | 47 |
| 31622.7767185595693… | n=2..388 | 387 |
| 55237.07504296764715433124781528617 | n=2..633 | 632 |
| 999982.6807693608… | n=1..899 | 899 |

The record number of values is 899 with c = 999982.6807693608… From which we could conjecture that for sufficiently large c, an arbitrary number of primes can be produced.

The 3rd formula is based on Fermat primes. Fermat primes are the primes of the form : $2^{2^n} + 1$. What is known so far is that for n = 1 to 4 we have the primes 2, 5, 17, 65537 and the pattern stops here apparently. Again, by modifying the premise we just have to write

$$F(n) = [r^{2^n} + {}^{1}/_{2}]$$

And to search for a suitable real value of r that will produce the first Fermat primes and the ones after (modified Fermat primes). If r = 2 then we get the usual Fermat numbers after which n > 4 fails to produce more primes. But if r = 1.0905077… then we get this series of primes:

2, 3, 5, 17, 257, 65537, 4294967311, 18446744193968636141,
340282371357715587431288126011714099603,
115792092256830257597513487698137234684227436353307972878385071833485576558709, …
(see appendix for more prime values and the precise value of r).

# Description of the algorithm and method

There are 3 steps

1) First we choose a starting value and exponent (preferably a rational fraction for technical reasons).
2) Use Monte-Carlo method with the Simulated Annealing, in plain english we keep only the values that show primes and ignore the rest. Once we have a series of 4-5 primes we are ready for the next step.
3) We use a formula for forward calculation and backward. The forward calculation is

   Forward : Next smallest prime to { $a(n)^e$ }.

   It is easy to find a probable prime up to thousands of digits. Maple has a limit of about 10000 digits on a Intel core i7 6700K, if I use PFGW I can get a probable prime of 1000000 digits in a matter of minutes.

   Backward : (to check if the formula works)

   Previous prime = solve for $x$ in $x^e - S(n+1)$. Where $S(n+1)$ is the next prime candidate. This is where e needs to be in rational form in order to solve easily in floating point to high precision using Newton-like methods.

## Conclusions

There are no proofs of all this, just empirical results. In practical terms, we have now a way to generate an arbitrary series of primes with (so far) a minimal growth function. The formulas are much smaller in growth rate than of the 2 historical results of Mills and Wright. Perhaps there is even a simpler formulation, I did not find anything simpler. In the appendix, the 50'th term of the sequence beginning with $10^{500} + 961$ is given, breaking the record of known series of primes in either a polynomial (46 values) or primes in arithmetic progression (26 values).

# Appendix

Value of a(0) for $a_{n+1} = a_n^{\frac{5}{4}}$.

a(0) = (2600 digits)

43.8046877158029348185966456256908949508103708713749518407406132875267041950609\
3664963337107816445739567943558075746653930014636424735710790652072233219699229\
3264945523871320676837539601640703027906901627234335816037601910027530786667742\
4699462324971177019386074328699571966485110591669754207733416768610281388237076\
1457089916472910767210053536764714815683021515348888782413463002196590898912153\
7807818768276618328571083354683602327177557419282082356357338008605137724420625\
6355959550708768409691378766548276776271408551872012425050013678129199364572781\
6813484896708864364878618010738012262326340551265300995594932556249327764325520\
8194746104359540189525699560278047888104384271808675933509955483655110570975535\
0785549077202528823776864085546060978256419523714505379764311017583547117340087\
4119475081308042239895607717353189481000696610499835674474037156885696125148581\
8456576179531740469553899959017477684394053845672138911212798487820145291551046\
7629440955925010855475725385626417446232721714485692213455756668913904742940022\
8509137974545537496866153147014899471237034869279337802008797927308327868037503\
8887734734110132599896243006668693634844327103981195681572642575883497205246107\
7249503005611993311614482718587460518670384095172358532168924929832717916742877\
5928083798545071506803190331910539166466084176778289275156307759698474412455449\
4235840260513685727576427182775664320370267582830712839903568563176995646262746\
1416682265277556094918638843788662398254929224611054007569004788919208610982286\
8199056832247214421586610021379123061307988471051671497283110156867428960680180\
5503966626145521465670913389257049550857812576171202023957487357339788493712676\
9426225000014887911004760565168355053802557466312278070529726060791106644597456\
7661803493305130350891168486525370221416972540649352435163491961811066911684870\
8294575729521396968670925609170964827383150968820007346616462008754667125912162\
9887393234505470764134731474624437303908697918037642878970154570903124414003971\
0357272367808690366664091433772132766429665349666355316481741675272445299039148\
2286893771949904585620265483705408090265640475448008548569022696419278537015273\
9052515666046613248284007441714998112109106552831729128554638499405603790836674\
4055864672832545083146225507711246656708938743897521843934561767824939322527151\
9446656406377385343422986304586480251980132774829846894886302934727512320855666\
7311057826570702055978747565264786935758838660918835790956716726895968781893980\
6197577089835533428250070076463861672738236689824429341266143008249630480529249\
5142244721354743123104963974936146293364440867656174513881761811216195861355409\
7226590982832455160532348892270853032526217953271203753936241844699444780852653\
926712756105661

50’th prime in the formula $a_{n+1} = a_n^{\frac{101}{100}}$ and $a_0 = 10^{500} + 961 + \epsilon, 0 < \epsilon < 0.5$

S(50)= (807 digits).

129729528971426122166658259081315435974871367309456840812055525509563976052536464197\
821936120784492089449745630948278142648656401758919926499683620493424145145363861773\
044716845814540511418289754542689191694327904116242782241131052138054549585683795895\
226460529926493834263717492409387560259409231253958370245042303023794648019244182073\
576593618946511947995963350548413770285593359081097306798650486731513585054871329096\
194202981055877907668708729761964242992640744211230936407662435884639367683685800000\
716124853576007781499789743771269181463159253173337794440878414346193538514506034277\
502087533266305538298562224619861085522581430515597209416207494298867400378422593043\
260350351208262898632520628116793338057678207643439460644660886621181985756002255888\
259043523402372168932260997906477619348535003398763

Primes generated by $\{c^n\}$, $n > 1$, will produce 387 consecutive primes. { } is the nearest integer.

c = 31622.776718559569341181978706142880921196050198706691685092604397231785036207 4\
0864754895287267371526503218374753273756995351706550976656637342668814334968666\
7101253409003404462577908490644122849185117910040141357956901545395665144254926\
0368733397963027260382237357834177624731483528295517610175141894026910108773091\
9102265276155731349632552235646978900096647191268683240475690792929946636829523\
4537889911695242127537567542061489388029765753953635000339676196085538110591110\
4585141090122658582123787181713420542296072149094895524374290222132677205020669\
4445695634733625075384860598871603780178145604934183538244723142272254031969991\
1200615639448998851471276974790874283370668013877448445827189836933336932748683\
7096031634918922648353422318092987818160592340775276127611835724115334227731181\
2382996597965951469528485537607311159096901622414506311547051905763888568046565\
1049166427400916041959872354321110387678917719194813843955193371342978053712862\
0453524374293366590642009861486573168440136568062290968471195962360800157373661\
7143373070165944776420832972019845288922599364819623823506773804799606015475838\
5408260449690988543205825633919510504581682546272086621725352722128830938671472\
3701049700231880709516101390100111272624859303900292400178495789671057286940063\
2939593946897204453503410798379799657214376636228419529879327743486265567548211\
3576902297196496895449727060689840512311047588147145276057014367506675470634720\
3613695309242488801761837728748016218351333849539796405079226399842615890863407\
8599788643344554273489131657070089852713689432765198710321770400528221171272704\
0206817976260112923629070254233895740091029484512516533519548133198851927540897\
2455463807670389342975063899265873368827496362107236544910631527268372333090949\
6599079804007385491895861312920245982318590485133202228187065916469134386843632\
6710041357132606931061476180276702785264896516822792415007247720747091646449625\
4232619886412790996612644285061938042714365958818635813884068102869316665794217\
8129882233115569739186048345925987950967419367699136707185345903901351463318825\
1244196261056665156495557988992317891312059314344240462240709984206444212348221\
8190871664878419254779689447073126235375450129126837015443848887673451873007831 7\
2181027964646479454712955925745105851307097731929480992845429539117442668639185\
0136062499136302371850993390404927892216296025369572878087609243226453797867139\
3469752835416017623841715414916352205819923096864111587928911940669604635015684\
7581294362465199335660467818710512634839702174274066131707754882547089414783999\
5175780316883752292172422813518529901439742936742350859348052219626008764448086\
8799805941436104497986393192697940430432240135606638442974285495155375736097977\
5028072301126799837273791743647186112604934117249206784451299546711413367340501\
9569662946589302409473910039374488077465339877358211311017648701799022641035886\
3484903089265972324361003850699420234653560932806297966170293267500259491821818\
1198369337121033424541848925401516832417018610965424291916655343744273160443306\
9511466370054097794456240488450136865019379631160402883720249277945984213446264\
5988545557033239241001377665704382870870778480345318077395858285768401082992990\
6805331632859321631952446031722731370128828153399370730874506904362805691466137\
7583386102156444256251269973458369382478957925560223823273029506595199057689583\
4494612415638418972310233727907333983390896921545293647542190449362825483734461\
9136791675945740812366000990257217886571673799707804765177555198609248755512965\
3027715120275248760325555105027812428659621914108249357373916718531633964814113\
5054938247521657216756907893392593294500494373382654321063600590692748010279682\
9489802744322026037605807751069630567877106566741622790692273950877423124270171\
5917401923410677957081175116161669874591767721421034406660019253785321121319776\
9175525445047016860655132021850376701748670306588930285165633137124920560978651\
3594698267301854434353463258169713980181998659759698957791522027302972704755488\
1930524920113346550223965908425088243203667733442482257635397243984316985885847\
8313649362480390679519963175568339592886446984299461893063001573950 9...

## Value of r for the Modified Fermat primes :

1. 090507732679143013260507434137871548876169139155253779637071740112482901941131794927614\
165046070727875947453049568456644638663863401452969304983696039368421966287954145103497 85\
455664518925506557480015624283354745541187800260237441434585818347722472782307786524703 01\
557607223687241626640266751272088743461716576567553383281405335074629322749596163614073 60\
489036963360216461159785107587006992362346239535666427851674987741633450324742973043111 56\
918186387564450658210917363379376495007263712050095542815271813730842051590178254450601 26\
752908663727993544583241249049592594152849795471622102290959653822850136216969080387991 96\
383765509946043160491654036964742400917664070054350091845594182058374371179120883957081 56\
209078863947050980206182296933121091613242790093184483373032553760846280097529146535026 78\
419385100092354406574326943906560907540707231181351944607974226054289879689424551236499 33\
269059229458421568428428015944756318589900397750269761852844216798990985594666449075870 71\
480852972785889053247716805445723363473127585208382346554817097179928110754768062838403 08\
066804082648385495810615341050875706050914159701065280480706360464815526935157288853036 30\
549938972715599267522628158998537812982857061583862766054101539552319187514985791390727 80\
139767752718563189026492873149585990085770880006556885215599134917493469603637739890197 37\

19405959699871790578939570721739263715332559979008810301240908542254896518354126910698473\
03993682210643911870348721611449094334837643545992733864032000150053382357189576849724304\
79394335286886902981957261182482082655523052080266649911319296914410341193709411800504767\
48354617111791391157262855664365330195306099000980663754432073067254803009304417954245510\
80052381553501227814722374446731494849665385984326298811570377471058095027488666777928067\
49282559951544564929643885617816410390414801853621166126904122739968733539945986005206956\
02391390456156785822775347992464020070450927826557980054516073143397824859022957041413936\
53297627437164784810222439767761820460207106318202527920152246071560668078195643604118308\
70602378959280839307869485035404878948607356189883649003438806376642735933157210850856805\
03506159529430925706576262450481876645606928209474175377071836512411835282091191296556778\
61072581015182894612172856050613393374195649787022411161276629365358921803831184186303787\
56578207686281678550906395477987278981782123483814461224106302196256888809443354614655931\
20883556394922774097730661062822409098069057549865237100948975335505690783704633949492808\
92879605033765196086164243915095640819228189066884662493174796111828932547748114051910717\
28593698556768965176931069943109211521660225588723269917817204936617401029568646551564000\
36551661428507983198000737013649470620201716629415463311634217579514585081448847531256578\
25714937231555792580624059121015209416549292439788842583514455246179201044493704583910868\
28830972654330281395445933979851895574225163914764337784816688407440642284922828049031522\
01618353388838173635532223004434950057652347341697358445971271679021336442636403592722174\
42375877255175408012442092930843055522454721678264477991680327705645021452815801955750023\
68386760695105940014524451622585268105939518510562367837875856568192370007002855068674279\
55506815221484917003793153080159210060716162049547569777826477385159377388028034148097703\
11536968521878819974014639448599782845579381209944477341673310914094632310493432676314786\
03847729990168248867135516819526306233740604368060524394007592069091131904709906212667687\
11025327697032309138725908179113445403295547528379573821310763160175025410458359533624152\
27347099777607748695939767186100079042186313608421185422966786497354923671602219894851764\
52621121094402997482563075776628177848989131288776952783081820635305570750177097902497128\
17820915639889334413242097535368258991208278875052331267293940341427764121856828064274269\
04307488474655053309132385768413627062332268645725260747542547126561064399074885662267080\
98171115376665516462361682889963839770000386660854204451896689392031192319488493110855357 5\
10266806327430894877402269963946041527164785009347817196388856921947720278211147199782688\

Modifies Fermat primes:

2, 3, 5, 17, 257, 65537, 4294967311, 18446744193968636141,
340282371357715587431288126011714099603,
115792092256830257597513487698137234684227436353307972878385071833485576558709, 134078086\
29214289698191742695172284384343525309923651738693154292959754483349760870492860683648451\
3268122749515235125357093412210435308279378947037706292 79, 17976933223763317017008802688 4\
50845059260725731988258572967021469957141212554036273433033527170307819182667579620360866\
32471077573785377675676242507349422625188943366906764576673141528562623002007942016657593\
90883867097480294726378517706587276267484299500204412355407141006057494372924502257219698\
352110801927, 323170128131645365815076255135382646669278924751793572195552962545141719652\
64801403285247181358559850418369614308791215550771771474749468184741737297690057465193509\
63945905691375542263397996520864769510659365204874119526029393967225952392638861315457919\
63108282494937753987330958790356284722814028831250184315941475044351819744600571075621618\
91978078908010514773644027856618091755427013769073948650806222690073031800766002745811262\
17470632677757168145724097916814902375242298300482723097426477384391005496248708782736397\
85136638055501928789129239806873288362297236918221684091632681230224992540958690935289708\
05311667, 10443893171662408345946053774430978239078973511588489688023962239907639155299 23\
35735754205395972104335511348468746703462926848540116993310568418686077261633016448147919\
18781556374086563714969978706718863915534446770744346798098445623101221088765942162306552\
14621897440640949470144818060504055261313257506351263985875277343435085681089703739446815\
48500161594788046362932860178367430582323227819014770251706417135952611655420880255303182\
64681512789697795902401025093747702623759547554737815750417887995127647950217987243038766\
77533890003260364907468957019674181713392537979549059382616660583120543334829823484146591\
47717501906341315211489961755543220463855382535950075700566535477413593697275201691979923\
38641568982946895521986824611585033511719052239621316679060986921373718762688888837745422\
83169367906229978113397931970222522988444158095482888380174612676484538739600000996356088\
45000911940693727397524975428611394844926595732425377211428922912540710834804687928863044\
33679420861353478682425113923228694593794164457617570785058604598399524014935395382862380\
98333161261112242273472110203395411585454959169383671956007018359097021603202728284219951\
102044249080006323598825172316922737869910163299754485996931556116349948454919395698791,
10907490458109667924271757791513592293983166654301156822150591558598181044442382010399174\
92919511601279745737961889572367787422187043894748772952513879077053306947607149615936661\
75112967878119912006205571981584069736977530603282183549613013123916787655598082588924149\
40891379369301309816262380113991728427333406116199557348489554551799814745807267443201348\
34707865566910864427699247698231429315759272436031068117019906918001830336446962391804107\
22985340156979073167257435430296713066681084447199038191801425151002140169249568182087355\
37216931832857484768370735100604211780047308319725147598436806135282678368780265876110257\
53988182574348514613527544858201515292421197926011282801578254357559923204391580750693084\
80293384143141225083638688402007089934075995374391110014446085562005376548363427140256578\
80908388428167581904401650025654343131473880484355345444054888207153825885498689967570953\
70985204324044035890167854980468740591429549828136949302322295576903613561409781403727774\
85357093055388039711174843200724936509378775007630503199092851376287125049320867683678108\
33352369655462724762257881165906356802040669563085578699320870099396037097662944064618881\
78983789996479926314834866382506382527233365361033178023818660797724861354693574334402892\
64201680765323750650601736988219125208588373661116687694701583189612511308498971041370105\
66940709481058594272986767876269978297997496028861741255819903164165855307471882021560488\
46419624409275268608448955329531354914194090887107690516242871900980302183047175728790288\
75725800629989849920195442450635678864794991294524049166732285923131275124915709630178011\

95007750822749047043672754341054932799971741515843977613033068821444197643295712208115235\
28096036461198603387546097428870083996323823306860603984960807123643767912944085743355903\
55150134317023096965637947743236134732979702590603391356175620665336125866408430315394702\
89835851661541255483770242672579803599405016279803061405113768381039638122682649588584888\
93930854366167295599946845670094424486542939255420837897352516799960987153631076804984709\
95101269541394419133543181784050825199638913349633307664893542738931263128355685853014270\
39519093232327765542481347851062836108789141693181881486942394111395782143969968340113780\
00491584853335395982836075034898639850515917970108430905151231503584598134574477059277161\
06628882069129509101520519702817443789175537317708652543820801180775690044296871487181349\
485685344318064027484651166891124430892293690879227355240747 3599